\documentclass[doublespacing,noinfoline]{imsart}

\RequirePackage[OT1]{fontenc}
\RequirePackage{amsthm,amsmath}
\RequirePackage{natbib}
\RequirePackage[colorlinks,citecolor=blue,urlcolor=blue]{hyperref}
\usepackage{graphicx}

\startlocaldefs
\numberwithin{equation}{section}
\theoremstyle{plain}

\newcommand{\X}{{{\cal X}}}

\newcommand{\Ob}{{{\cal O}}}
\newcommand{\bon}{{\bar {\cal O}_n}}

\newcommand{\Eet}{{\rm E}_{*}}
\newcommand{\Ee}{{\rm E}}

\newcommand{\var}{{\rm var}_{*}}

\newcommand{\dpg}[2]{\frac{\partial #1}{\partial #2}}
\newcommand{\ddpg}[2]{\frac{\partial^2 #1}{\partial #2^2}}

\newcommand{\PP}{{{\rm P}}}
\newcommand{\QQ}{{{\rm Q}}}

\newcommand{\bn}{\hat \beta_n}
\newcommand{\gn}{\hat \gamma_n}
\newtheorem{Definition}{Definition}
\newtheorem{Example}{Example}

\newcommand{\dd}{{\rm d}}
\newcommand{\AIC}{{\rm AIC}}

\newcommand{\KL}{{\rm KL}}
\newcommand{\EKL}{{\rm EKL}}

\newcommand{\CE}{{\rm CE}}
\newcommand{\ECE}{{\rm ECE}}
\newcommand{\CVCE}{{\rm CVCE}}

\newcommand{\HH}{{\rm H}}
\newcommand{\EH}{{\rm EH}}
\newcommand{\I}{{\rm I}}

\newcommand{\nn}{n^{-1}}
\newcommand{\si}{\sum_{i=1}^n}

\newcommand{\tendp}{\longrightarrow \hspace{-0.50cm} ^p \hspace{0.50cm}}

\endlocaldefs

\begin{document}

\begin{frontmatter}
\title{Information Theory and Statistics: an overview}
\runtitle{Information Theory and Statistics}

\begin{aug}
\author{\fnms{Daniel} \snm{Commenges}\ead[label=e1]{daniel.commenges@isped.u-bordeaux2.fr}}

\address{Daniel Commenges\\ Epidemiology and Biostatistics Research Center, INSERM\\ Bordeaux University\\146 rue L\'eo
Saignat, Bordeaux, 33076, France\\
\printead{e1}}

\runauthor{D. Commenges}

\affiliation{INSERM and University Victor Segalen Bordeaux 2}

\end{aug}

\begin{abstract}
We give an overview of the role of information theory in statistics, and particularly in biostatistics.
We recall the basic quantities in information theory; entropy, cross-entropy, conditional entropy, mutual information and  Kullback-Leibler risk. Then we examine the role of information theory in estimation theory, where the log-klikelihood can be identified as being an estimator of a cross-entropy. Then the basic quantities are extended to estimators, leading to criteria for estimator selection, such as Akaike criterion and its extensions. Finally we investigate the use of these concepts in Bayesian theory; the cross-entropy of the predictive distribution can be used for model selection; a cross-validation estimator of this cross-entropy is found to be equivalent to the pseudo-Bayes factor.

\end{abstract}


\end{frontmatter}

\section{Introduction}
 \citet{shanon1949mathematical} introduced the concept of {\em information} in communication theory. One of he key concept in information theory is that of {\em entropy}; it first emerged in thermodynamics in the 19$^{th}$ century and then in statistical mechanics where it can be viewed as a measure of {\em disorder}. Shannon was searching for a measure of information for the observation of a random variable $X$ taking $m$ different values $x_j$ and having a distribution $f$ such that $f(x_j)=\PP (X=x_j)=p_j$.  He found that  entropy was the only function satisfying three natural properties: i) $H(X)$ is positive or null; ii) for a given $m$, the uniform distribution maximizes $H(X)$; iii) entropy has an additive property of successive information. The last property says that if one gets an information in two successive stages, this is equal to the sum of the quantities of information obtained at each stage. For instance one may learn the outcome of a die roll by first learning whether it is even or odd, then the exact number. This is formula (\ref{additivity}) given in section \ref{condentropy}. Entropy is defined as :
\begin{equation} \label{dentropy} \HH(X)=\sum_{j=1}^m p_j\log \frac{1}{p_j}.\end{equation}
The quantity $\log \frac{1}{p_j}$ measures the ``surprise'' associated with the realization of the variable at the value $x_j$, and $\HH(X)$ is the expectation of $\log \frac{1}{f(X)}$.  $\HH(X)$ can be interpreted as measuring the quantity of information brought by $X$.
Taking base-2 logs, this defines the information unit as the information brought by a uniform binary variable with probabilities $p_1=p_2=0.5$.
 \citet{khinchin1956basic} studied the mathematical foundations of the theory of information; a more recent book is that of  \citet{cov91}, and an accessible introduction can be found in \citet{applebaum1996probability}.

The use of information theory was introduced in statistics by \citet{Kul51} and developed by Kullback in his book \citep{Kul59}. In statistics, entropy will be interpreted as a measure of uncertainty or of {\em risk}. Although entropy is defined by the same formula in physics, communication theory and statistics, the interpretation is different in the three domains, being respectively a measure of disorder, of information, of uncertainty. Moreover in communication theory the focus is on entropy for discrete variables while in statistics there is often a strong interest in continuous variables, and the interpretation of entropy for continuous variables, called ``differential entropy'', is not obvious. The {\em relative entropy} or Kullback-Leibler risk, as well as the related concept of {\em cross-entropy},  play a central role in statistics. The aim of this paper is to give an overview of applications of information theory in statistics, starting with the the basic definitions. More technical reviews can be found in \citet{ebrahimi2010information} and \citet{clarke2014statistical}.

In section 2 we review the basic concepts of information theory with a focus on continuous variables, and discuss the meaning of these concepts in probability and statistics. In section 3 we see that the maximum likelihood method can be grounded on the Kullback-Leibler risk, or equivalently, on cross-entropy. In section 4 we show that this is also the case for some of the most used criteria for estimator selection, such AIC and TIC. In section 5 we look at the use of information theory in the Bayesian approach, with two applications: measuring the gain of information brought by the observations, and model selection. Section 8 concludes.

\section{Basic definitions and their interpretation}

Conventional quantities in information theory are the entropy, the Kullback-Leibler divergence, and the cross-entropy. We shall focus on continuous variables but most of the formulas are also valid for discrete variables.

\subsection{Entropy}
For a continuous variable with probability density function $f$, one may define the so-called ``differential entropy'':
\begin{Definition}[Differential entropy]
$$\HH(X)=\int f(x)\log \frac{1}{f(x)} \dd x.$$
\end{Definition}
It can be called simply ``entropy''.
Note that this is a function of $f$ so that it can be written $\HH(f)$, and this is more relevant when we consider different probability laws.
For continuous variables, base-2 logarithms do not make particular sense and thus natural logarithms can be taken. As in the discrete case, $\HH(X)$ can be viewed as the expectation of the quantity $\log \frac{1}{f(X)}$, that can be called a {\em loss}. In decision theory a {\em risk} is the expectation of a {\em loss}, so that $\HH(X)$ can be called a risk. This differential entropy has the additive property (\ref{additivity}) described in Section \ref{condentropy}, but its interpretation as a quantity of information is problematic for several reasons, the main being that it can be negative. It may also be problematic as a measure of uncertainty.

It is instructive to look first at the distributions of maximum entropy, that is, the distributions which have the maximum entropy among all distributions on a given bounded support or on an unbounded support with given variance. For the three most important cases, bounded interval $[a,b]$, $(-\infty,+\infty)$ and $(0,+\infty)$, these are respectively the uniform, normal and exponential distributions. Table $\ref{maxentdis}$ shows the entropy of these distributions.
Noting that the variance of $U(a,b)$ is $\sigma^2=1/12 (b-a)^2$, the entropy is $1/2 \log 12+ \log \sigma$.
Noting that the variance of the exponential is $\sigma^2=\lambda^{-2}$, the entropy is $1+\log \sigma$.
For the normal distribution the entropy can be written $1/2 \log (2\pi e)+\log \sigma$. That is, the entropy of these maximum entropy distributions can be written $\log \sigma$ plus a constant. For some other unimodal distributions we have also this relation; for instance the Laplace distribution has entropy $1+ 1/2 \log 2 + \log \sigma$. Thus, the entropy often appears to be equivalent to the variance in the log-scale.

However, this is not the case in general. For multimodal distributions there can be large differences in the assessment of uncertainty by using  variance or entropy; for instance we can easily find two distributions $f_1$ and $f_2$ such as $f_1$ has a smaller variance but a larger entropy than $f_2$. In that case the question arises to know which of the two indicators is the most relevant. While it may depend on the problem, one argument against entropy is that it does not depend on the order of the values, while in general the values taken by a continuous variable are ordered. A strength of the entropy is that it can describe uncertainty for a non-ordered categorical variable where the variance has no meaning. On the contrary, when the order does have a meaning the variance may be more informative than the entropy.

\begin{table}
\caption{Maximum entropy distributions for three cases of support: in the three cases, the maximum entropy is a constant plus the logarithm of the standard deviation ($\sigma$). }\label{maxentdis}
\begin{center}
\begin{tabular}{rrrr}
  \hline
Support & Distribution & Entropy & $C+ \log \sigma$ \\
  \hline
 Bounded interval $(a,b)$ & Uniform & $\log (b-a)$ & $1/2 \log 12+ \log \sigma$ \\
  $(-\infty,+\infty)$ with given variance  & Normal  & $1/2 \log (2\pi e\sigma^2)$ & $1/2 \log (2\pi e)+\log \sigma$  \\
  $(0,+\infty)$ with given variance & Exponential & $H=1-\log \lambda$  & $1+\log \sigma$  \\
     \hline
\end{tabular}
\end{center}
\end{table}

\begin{Example}[Will there be a devastating storm ?]
A storm is forecast: let $X$ be the speed of the wind in km/h. Compare two distributions: i) $f_1$, a normal distribution with expectation $100$ and standard deviation $10$; $f_2$, a mixture distribution, with weights one half, of two normal distributions with standard deviation $3$ and expectations $50$ and $200$. It is clear that $f_1$ has a smaller variance but a larger entropy than $f_2$, because $f_2$ is quite concentrated on two particular values ($50$ and $200$). However from a practical point of view there is much more uncertainty with $f_2$ (we don't know whether there will be ordinary wind or a devastating storm) than with $f_1$ (there will be a strong but not devastating storm).

If we consider that a devastating storm is defined by a wind speed larger than $150$ km/h, one way to reconcile the points of views is to say that the important variable for taking a decision is not $X$ but $Y=I_{X>150}$. It is clear that both entropy and variance of $Y$ are larger under $f_2$ than under $f_1$.
\end{Example}

\subsection{Cross-entropy and relative entropy or Kullback-Leibler risk}
When there are two distributions $g$ and $f$ for a random variable $X$, we can define the cross-entropy of  $g$ relative to $f$.
For a discrete variable taking values $x_j, j=1,\ldots,m$, and denoting $p_j=f(x_j)$ and $q_j=g(x_j)$,  this is:
\begin{equation*}\label{crossdentropy} \CE(g|f)=\sum_{j=1}^m p_j\log \frac{1}{q_j}.\end{equation*}
It can be viewed as the expectation of a ``surprise indicator'' $\log \frac{1}{g(X)}$ under the distribution $f$. It can be decomposed in:
\begin{equation*} \CE(g)=\sum_{j=1}^m p_j\log \frac{p_j}{q_j}+\sum_{j=1}^m p_j\log\frac{1}{p_j}.\end{equation*}
Here $\sum_{j=1}^m p_j\log \frac{p_j}{q_j}$ is the {\em relative entropy} or {\em Kullback-Leibler risk} of $g$ relative to $f$, and is denoted $\KL(g|f)$. Here, the interpretation in terms of uncertainty, or ``risk'' is more relevant than that in terms of information. We can tell that, if the observations come from $f$,  the risk associated to $g$ is:
\begin{equation}\label{KL+H} \CE(g|f)=\KL(g|f)+\HH(f),\end{equation}
that is, this is the sum of the Kullback-Leibler risk (of using $g$ in place of $f$), and the entropy of $f$ (the risk already associated to $f$).

This extends to the differential cross-entropy (for continuous variables).

\begin{Definition}[Cross-entropy]
 The cross-entropy of $g$ relative to $f$ is:
\begin{equation}\label{cross-entropy} \CE(g|f)=\int f(x)\log \frac{1}{g(x)} \dd x.\end{equation}
\end{Definition}

\begin{Definition}[Kullback-Leibler risk]
The Kullback-Leibler risk of $g$ relative to $f$ is:
\begin{equation}\label{Kullback-Leibler} \KL(g|f)=\int f(x)\log \frac{f(x)}{g(x)} \dd x.\end{equation}
\end{Definition}
The decomposition \ref{KL+H} is still valid for these definitions.

Different notations and names have been used for the Kullback-Leibler risk; it is often noted $D(f||g)$ where ``D'' stands for  ``divergence''. It is sometimes called ``distance''. Indeed it is a measure of how far $g$ is from $f$ and it has some of the properties of a distance: for all $f, g$ we have  $\KL(g|f) \ge 0$; we have $\KL(f|f)=0$, and conversely if  $\KL(g|f)=0$ then $f=g$ nearly everywhere. However it is not a distance because it is not symmetric and does not satisfy the triangular inequality. A symmetrized version satisfies the distance properties but we will stay with the dissymmetric version which is well adapted to the statistical problem. In fact, this dissymmetry feature is important: $f$ specifies the reference probability measure which is used for taking the expectation of $\log \frac{f(X)}{g(X)}$, and in statistics the reference probability measure will be taken as the ``true'' one.

The Kullback-Leibler risk is invariant for affine transformations of the random variable: if one defines $Y=aX+b$ , we have $\KL(g_Y|f_Y)=\KL(g_X|f_X)$, where  $g_Y$ and  $f_Y$ are the transformed density for $Y$.

\subsection{Conditional and expected conditional entropy}\label{condentropy}
The conditional entropy is defined by the same formula as the (differential) entropy, but using the conditional density $f_{Y|X}$:
$$ \HH^X(Y)=\int f_{Y|X}(y|X)\log \frac{1}{f_{Y|X}(y|X)} \dd y.$$
With this definition the conditional entropy $\HH^X(Y)$ is a random variable, and it can be smaller or larger than the entropy $\HH(X)$.
This is not conventional however: what is generally called conditional entropy is the expectation of this quantity, here called {\em expected} conditional entropy:
 $$ \EH^X(Y)=E[\HH^X(Y)]=\int \int f_{Y|X}(y|X=x)\log \frac{1}{f_{Y|X}(y|X=x)} \dd y ~~f_X(x) \dd x.$$

The expected conditional entropy is always lower than the entropy as we shall see in section \ref{mutualinformation}. We show an example with a discrete variable.

It is interesting to make the distinction between conditional entropy and {\em expected} conditional entropy because the conditional entropy may be useful in applications; \citet{alonso2007surrogate} did make this distinction. Although the {\em expected} conditional entropy is always lower than the entropy, this is not the case for the conditional entropy.
 \begin{Example}[Will your plane crash ?]
 Suppose you are in a plane and $Y=1$ if the plane crashes before landing and $Y=0$ otherwise. Variable $X$ represents the state of the engine with $X=1$ if an engine is on fire and $X=0$ otherwise. Assume that $\PP(X=1)=0.001$, $\PP(Y=0|X=0)=1$ while  $\PP(Y=0|X=1)=0.5$. The entropy of $Y$ is very small because there is almost certainty that the plane will not crash: $\PP(Y=0)=1\times 0.999+ 0.5\times 0.001=0.9995$. We find $H(Y)=\PP(Y=0)\log \frac{1}{\PP (Y=0)}+\PP(Y=1)\log \frac{1}{\PP (Y=1)}=0.0043$

 The conditional entropy given that $X=1$ (engine in flame) reaches the maximum for a binary variable: $\HH(Y)=1$. However the expected conditional entropy is lower than the original entropy because the probability that $X=1$ is very small: $\EH^X(Y)=0 \times \PP(X=0) + 1 \times 0.001=0.001$, which is indeed smaller than $0.0043$. In practice, if you really learn that $X=1$ you would like to use the conditional entropy which better describes the uncertainty of the flight at this moment.
 \end{Example}

A way to quantify the information gain brought by $X$ is to measure the distance between the marginal and the conditional distribution. If we take the conditional distribution  as reference, the information gain can be defined by
$$ {\rm IG}(X\rightarrow Y)=\KL[f_Y|f_{Y|X}(.|X)].$$
In contrast with the entropy change $\HH(Y)-\HH^X(Y)$, the information gain is always positive or null. However, as displayed in Equation (\ref{expInfGain}) below, their expectations are equal.

\subsection{Joint entropy and mutual information}\label{mutualinformation}
The joint entropy of X and Y is :
\begin{equation}  \HH(X,Y)=\int \int f_{Y,X}(y,x)\log \frac{1}{f_{Y,X}(y,x)} \dd y \dd x.\end{equation}
We have the following additive property:
 \begin{equation}\HH(X,Y)=\HH(X)+\EH^X(Y).\label{additivity}\end{equation}
 This is an essential property that we ask for a measure of information: the information in observing $X$ and $Y$ is the same as that of observing $X$ and then $Y$ given $X$. This makes also sense if the interpretation is uncertainty.

\begin{Definition}[Mutual information]
 The mutual information is the Kullback-Leibler risk of a distribution defined by the product of the marginals relative to the joint distribution:
 \begin{equation}\label{mutual}  \I(X;Y)=\int \int f_{Y,X}(y,x)\log \frac{f_{Y,X}(y,x)}{f_X(x)f_Y(y)} \dd y \dd x.\end{equation}
 \end{Definition}

 If $X$ and $Y$ are independent, the mutual information is null; thus, $I(X,Y)$ can be considered as a measure of dependence.
 We have the following relation:
 \begin{equation}\I(X;Y)=\HH(Y)-\EH^X(Y).\end{equation}
 It follows that $\EH^X(Y)=\HH(Y)-\I(X,Y)$. Thus, since $\I(X;Y)\ge 0$, the expected conditional entropy of $Y$ is lower than its entropy.

 Mutual information can also be expressed as the expectation of the information gain (the Kullback-Leibler risk of $f_Y$ relative to $f_{Y|X}$); the mutual information can be written $\I(X,Y)=\int \int f_{Y,X}(y,x)\log \frac{f_{Y|X}(y|x)}{f_Y(y)} \dd y \dd x.$ so that we have:
 \begin{equation}\label{expInfGain} \I(X;Y)=\Ee \{ \KL[f_Y|f_{Y|X}(.|X)]\}=\Ee [\HH(Y)-\HH^X(Y)]=H(Y)-\EH^X(Y).\end{equation}
 Thus, both the expected entropy change and the expected information gain are equal to the mutual information.
Mutual information was used by \citet{alonso2007surrogate} to quantify the information given by a surrogate marker $S$  on the true clinical endpoint $T$. They propose the measure $R^2_h=1-e^{-I(S;T)}$ which takes values between $0$ and $1$. For normal distributions this measure amounts to the ratio of the explained variance and the marginal variance of $T$.

We can now define the conditional mutual information:
\begin{Definition}[Conditional mutual information]
The conditional mutual information is:
 \begin{equation}\label{condmutual}  \I(X;Y|Z)=\int \int f_{Y,X|Z}(y,x|Z)\log \frac{f_{Y,X|Z}(y,x|Z)}{f_{X|Z}(x)f_{Y|Z}(y|Z)} \dd y \dd x.\end{equation}
\end{Definition}

We have the additive property:
\begin{equation} \label{additiveinfo} \I(X,Z; Y)=\I(Y;X)+\I(Y;Z|X).\end{equation}
This can be generalized in the so-called ``chain-rule for information'' \cite{cov91}.
Thus, there is additivity of information gains: the information brought by $(X,Z)$ on $Y$ is the sum of the information brought by $X$ and the information brought by $Z$ given $X$.

\begin{Example}[Mutual information in a normal regression model]
Assume we have:
$Y=\beta_0^*+\beta_1^*X+\beta_2^*Z+\varepsilon$,
with $X$ and $Z$ and $\varepsilon$ normal with variances $\sigma_X^2$, $\sigma_Z^2$ and $\sigma_{\varepsilon}^2$;
$\varepsilon$ independent of $X$ and $Z$ and ${\rm cov}(X,Z)=\sigma_{XZ}$.
The marginal and conditional distributions on $X$ and $X,Z$ are normal with variances respectively:
$$\var(Y)=\sigma^2_{\varepsilon}+\beta_1^2\sigma_X^2+\beta_2^2\sigma_Z^2-2\beta_1\beta_2\sigma_{XZ}$$
$$\var^X(Y)=\sigma^2_{\varepsilon}+\beta_2^2\sigma_Z^2$$
$$\var^{X,Z}(Y)=\sigma^2_{\varepsilon}$$
Since for a normal distribution the entropies are $\log \sigma +c$, the predictability gains are
\begin{eqnarray*}\HH(Y)-\EH^X(Y)&=&\frac{1}{2}\log (\var(Y))-\frac{1}{2}\log (\var^X(Y))\\
&=&\frac{1}{2} \log \frac{\var(Y)}{\var^X(Y)}\end{eqnarray*}
\begin{eqnarray*}\HH(Y)-\EH^{X,Z}(Y)&=&\frac{1}{2}\log (\var(Y))-\frac{1}{2}\log (\var^{X,Z}(Y))\\
&=&\frac{1}{2} \log \frac{\var(Y)}{\var^{X,Z}(Y)}\end{eqnarray*}
For normal distributions the information gains are the logarithm of the factor of reduction of standard deviations, and it is easy to verify the additive property (\ref{additiveinfo}).
\end{Example}

\subsection{Extensions: definitions for probability laws and sigma-fields}\label{extensions}
The Kullback-Leibler risk can be defined for general probability laws as in \citet{Kul59} and \cite{haussler1997mutual}. For general distributions $\PP$ and $\QQ$, the Kullback-Leibler risk of $\QQ$ relative to $\PP$ is
defined by: $\KL(\QQ|\PP)= \int \dd \PP \log \frac{\dd \PP}{\dd \QQ}$, where $\frac{\dd \PP}{\dd \QQ}$ is the Radon-Nikodym derivative of $\PP$ relative to $\QQ$.

The Kullback-Leibler risk has also been studied among other metrics for probabilities in \citet{gibbs2002choosing}. It has been used as a metric for some asymptotic results in statistics \citep{hall1987kullback,clarke1990information}.

Since Radon-Nikodym derivatives are defined for sigma-fields, we can make explicit that the Kullback-Leibler risk also depends of the considered sigma-field. So, the Kullback-Leibler risk of $\QQ$ relative to $\PP$ on sigma-field $\X$ can be noted:
$$\KL(\QQ|\PP;\X)= \int \dd \PP \log \frac{\dd \PP}{\dd \QQ}_{|\X}.$$
This notation and concept has been in used by \citet{liquet2011choice} to define a so-called restricted AIC, described in section \ref{ext-crit}.

\section{Information theory and estimation}
A statistical model for $X$ is a family of distributions indexed by a parameter:
$(f^{\theta})_{\theta\in \Theta}$. For parametric models $\Theta \subset \Re^d$.  We assume that there is a {\em true} distribution, $f^*$, which is unknown. The model is {\em well specified} if there exists $\theta_*$ such that $f^{\theta_*}=f^*$ ($f^*$ belongs to the model); it is  {\em misspecified} otherwise. The aim of inference is to choose a distribution as close as possible to the true distribution. The Kullback-Leibler risk is a good candidate to quantify how close a distribution $f^{\theta}$ is from the true distribution $f^*$. Here the dissymmetric property of the Kullback-Leibler risk is not a problem here because the situation is truly dissymmetric: we wish to quantify the divergence of $f^{\theta}$ relative to $f^*$. So we wish to minimize
$$\KL(f^{\theta}|f^*)= \Eet [\log \frac{f^*(X)}{f^{\theta}(X)}],$$
where $\Eet$ means that the expectation is taken relative to $f^*$.
Note that if the model is well specified the risk is minimized by $f^{\theta_*}$, and the minimized risk is zero. If the model is misspecified, there is a $f^{\theta_0}$ which minimizes the risk: $\KL(f^{\theta_0}|f^*)$ can be called the ``misspecification risk'' of the model.

The problem, of course, is  that the risk depends on the unknown $f^*$. Using the relation \ref{KL+H} and since the entropy of $f^*$ does not depend on $\theta$, we can simply minimize the cross-entropy, $\CE(f^{\theta}|f^*)= \Eet [\log \frac{1}{f^{\theta}(X)}]$. Thus, the risk $\CE(f^{\theta}|f^*)$ depends on $f^*$ only by the expectation operator. This is really a ``risk'', defined as the expectation of a loss function: the loss here is $\log \frac{1}{f^{\theta}(X)}$ which does not depend itself of $f^*$. If a sample of independently identically distributed variables $\bar \Ob_n=(X_1,\ldots,X_n)$  is available, the cross-entropy can be estimated by
$-\nn \sum_{i=1}^n \log f^{\theta}(X_i)$.
 We recognize it as the opposite of the normalized loglikelihood. By the law of large numbers this quantity converges toward $\CE(f^{\theta}|f^*)$:
 $$-\nn L_n(f^{\theta}) \tendp \CE(f^{\theta}|f^*),$$
 where $L_n(f^{\theta})$ is the loglikelihood.
 Thus, the maximum likelihood estimator (MLE) minimizes the natural estimator of the cross-entropy. If the model is well specified, the variance of the asymptotic distribution of the  MLE is given by the Fisher information matrix at $\theta_*$. It is easy to see that we have $\I(\theta_*)=n \ddpg{\KL(f^{\theta})}{\theta}|_{\theta_*}$: the individual Fisher information matrix is thus the curvature of the Kullback-Leibler risk (or equivalently of the cross-entropy).

\section{Estimator selection: normalized AIC and TIC}

Another important problem in statistics is that of model choice, or estimator selection.
The concepts of entropy and Kullback-Leibler risk have been adapted to describe the risk of an estimator of the true density function $f^*$ \citep{commenges2008estimating}. The need of an extension comes from the fact that the Kullback-Leibler risk is defined for a fixed distribution, while an estimator is random. Denoting  $g^{\bn}$ such an estimator, the expected Kullback-Leibler risk can be defined as:
 $$\EKL(g^{\bn})=\Eet [\log \frac{f^*(X)}{g^{\bn}(X)}].$$

 In the notation $\EKL(g^{\bn})$ it is implicit that the expected risk is wrt $f^*$. Similarly we define the expected cross-entropy
 $$\ECE(g^{\bn})=\Eet [\log \frac{1}{g^{\bn}(X)}].$$
  As the cross-entropy, the expected cross-entropy can be additively decomposed as
  $$\ECE(g^{\bn})=\EKL(g^{\bn})+\HH(f^*),
  $$ with the interpretation that the risk with $g^{\bn}$ is the risk incurred in using $g^{\bn}$ in place of $f^*$ plus the risk using $f^*$.

If the model is well specified, there exist a value $\beta_*$ such that $g^{\beta_*}=f^*$; in that case, $\EKL(g^{\bn})$ is the ``statistical risk''. It converges toward zero since $\bn \rightarrow 0$. Since
$\KL(g^{\beta_*}|g^{\beta_*})=0$ and  $\dpg{\KL}{\beta}|{\beta_*}=0$ we have for a fixed sequence $\{\bn\}$:
$$\KL(g^{\bn}|g^{\beta_*})=1/2 (\bn-\beta_*)^\top \ddpg{KL(g^{\beta}|g^{\beta_*})}{\beta}|_{\beta_*}(\bn-\beta_*)+o(||\bn-\beta_*||^2).$$
Taking expectation and since the asymptotic variance of $(\bn-\beta_*)$ is $I^{-1}$ and that $\ddpg{KL(g^{\beta}|g^{\beta_*})}{\beta}|_{\beta_*}=I/n$, we find:
$$\EKL(g^{\bn})=p/2n+ o(\nn).$$

\begin{figure}
\begin{center}
\includegraphics[scale=0.45]{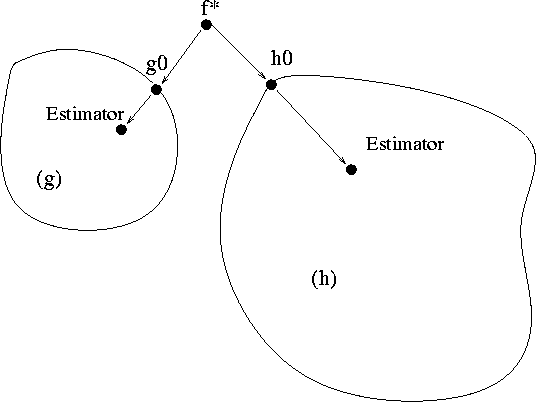}
\caption{Risk of two models $(g)$ and $(h)$ relative to the true distribution $f_*$. The risk of the estimator is the sum of the misspecification risk and the statistical risk.  }\label{models}
\end{center}
\end{figure}

If the model is misspecified, the estimator converges toward $g^{\beta_0}$ which is different from $f^*$. By definition, the MLE converges toward $g^{\beta_0}$ which has the smaller Kullback-Leibler risk in the
model.
The expected Kullback-Leibler risk of the MLE can be decomposed as:
$$\EKL(g^{\bn})=\Eet [\log \frac{g^{\beta_0}}{g^{\bn}}]+\KL(g^{\beta_0}|f^*),$$
where $\Eet [\log \frac{g^{\beta_0}}{g^{\bn}}]$ is the statistical risk in the misspecified case, and $\KL(g^{\beta_0}|f^*)$ is the ``misspecification risk''. Figure \ref{models} illustrates the risks for two non-overlapping models.

The statistical risk still tends to zero since $\bn \rightarrow \beta_0$. However there is no guarantee that  $p/2n$ is a good estimate of it. The misspecification risk is fixed. It is not possible to have a good estimator of $\KL(g^{\beta_0}|f^*)$ but we have an estimator of $\CE(g^{\beta_0}|f^*)$ by: $-\nn\sum_{i=1}^n \log g^{\bn}(X_i)$. This estimator however is biased downward because the MLE minimizes $-\nn\sum_{i=1}^n \log g^{\beta}(X_i)$. Computations (either direct or by cross-validation) show that the bias is $p/2n +o(\nn)$ \citep {Lin86}. If we estimate the statistical risk by $p/2n$ (as in the case of well specified models) we end up with the estimator:
$$-\nn\sum_{i=1}^n \log g^{\bn}(X_i)+p/n=\frac{1}{2n} \AIC.$$
Thus, the normalized Akaike criterion (AIC) \citep{Aka73} is an estimator of the expected cross-entropy ECE. \citet{takeuchi1976distributions} has proposed another criterion (TIC) which does not assume that the statistical risk can be estimated by $p/2n$.
Other criteria are available when the estimators are not MLE but general M-estimators, in particular the generalized information criterion (GIC) \citep{konishi1996generalised}. All these estimates of ECE are equal to minus the normalized loglikelihood plus a correction term.

Of course, estimation of the expected cross-entropy is essentially useful for choosing an estimator when several are available. The difference of the expected cross-entropy of two estimators $g^{\bn}$ and $h^{\hat \gamma_n}$ is equal to the difference of their expected Kullback-Leibler risks, so that the difference of estimators of ECE also estimates the difference of expected Kullback-Leibler risks:
$$[{\widehat \ECE}(g^{\bn})- {\widehat \ECE}(h^{\gn})]- [\EKL(g^{\bn}) - \EKL(h^{\gn})] \tendp 0.$$
Note that the difference of expected Kullback-Leibler risks is not random. However it is not fixed since it depends on $n$; we have $\EKL(g^{\bn}) - \EKL(h^{\gn}) \rightarrow \KL(g^{\beta_0}|f^*) - \KL(h^{\gamma_0}|f^*)$.

Moreover, if $g^{\beta_0}\ne h^{\gamma_0}$, $\sqrt{n} [{\widehat ECE}(g^{\bn})- {\widehat \ECE}(h^{\gn})]$ has a normal asymptotic distribution with a variance that can be easily estimated, leading to the construction of so-called ``tracking intervals'' \citep{commenges2008estimating} for the difference of expected Kullback-Leibler risks. The condition $g^{\beta_0}\ne h^{\gamma_0}$ may not hold if the models are nested: $(g) \subset (h)$. The hypothesis that the best distribution is in the small model ($g^{\beta_0} = h^{\gamma_0}$) is the conventional null hypothesis; it can be tested by a likelihood ratio statistic which follows a chi-squared distribution if the models are well specified. The distribution of the likelihood ratio statistics has been studied in the general setting by \citet{vuong1989likelihood}; these results are relevant here because the asymptotic distributions of the estimators of ECE is driven by that of the likelihood ratio statistics.

\subsection{Extensions: restricted and prognostic criteria}\label{ext-crit}
\subsubsection{Restricted AIC}
\citet{liquet2011choice} tackled the problem of choosing between two estimators based on different observations.  As one of their examples, they compared two estimators of survival distributions: one based on the survival observation, the other one including the additional information on disease events. Both estimators are judged on the sigma-field generated by observation of survival, noted $\Ob'$. The estimator which includes disease information uses more information, thus is based on a larger sigma-field than the direct survival estimator, and it will be judged only on the survival sigma-field; for defining this, it is useful to specify the sigma-field on which the Kullback-Leibler risk is defined, as suggested in Section \ref{extensions}. The difference of expected Kullback-Leibler risks on this restricted sigma-field was assessed by a criterion called $D_{\rm RAIC}$.
\subsubsection{Choice of pronostic estimators}
Prognostic estimators also are judged on a sigma-field which may be different from that generated by observations on which they are based. The problem is to give an estimator of the distribution of the time to an event posterior to a time $t$, based on observations prior to time $t$. Here again the sigma-field on which the estimator will be judged is not the same as the sigma-field on which the estimator has be been built; moreover there is the conditioning on the fact that the event has not occurred before $t$. \citet{commenges2012choice} defined  an ``expected prognosis cross-entropy'', or equivalently up to a constant, an ``expected prognosis Kullback-Leibler risk'' and proposed an estimator of this quantity based on cross-validation (called EPOCE).

\section{Information theory and the Bayesian approach}
In the Bayesian approach we may use the concepts of information theory to measure the information (or the decrease of uncertainty)  on the parameters or on the predictive distribution.
\subsection{Measuring information on the parameters}
\citet{lindley1956measure} was the first to use information theory to quantify information brought by an experiment in a Bayesian context. In the Bayesian approach, the observation of a sample $\bar \Ob_n$ brings information on the parameters $\theta$ which are treated as random variables. It seems natural to compare the entropy of the prior distribution $\pi(\theta)$  and of the posterior distribution $p(\theta|\bar \Ob_n)$.  Lindley's starting point was to inverse the concept of entropy by considering that the opposite of the entropy was a quantity of information; thus,  $\int \pi(\theta) \log \pi(\theta) \dd \theta$ and  $\int p(\theta|\bar \Ob_n) \log p(\theta|\bar \Ob_n) \dd \theta$ was considered as the amount of information on $\theta$ for the prior and the posterior, respectively.
We may then measure the information brought by the observations on the parameters by the difference of amounts of information between the posterior and the prior, or equivalently by the difference of entropy between the prior and the posterior distribution:
$$\Delta \HH(\theta,\bar \Ob_n)= \HH[\pi(\theta)]- \HH[p(\theta|\bar \Ob_n)].$$
Note that $p(\theta|\bar \Ob_n)$ is a conditional distribution, so that $\HH(p(\theta|\bar \Ob_n)$ is a conditional entropy. As we have seen, it is not guaranteed that $\Delta H(\theta,\bar \Ob_n)$ is positive. We may look at the expectation of $\Delta \HH(\theta,\bar \Ob_n)$. In virtue of \ref{mutual} this expectation is always positive and is equal to the mutual information between $\theta$ and ${\bar \Ob_n}$. However for a Bayesian this is not completely reassuring; moreover when we have observed $\bar \Ob_n$, we know the conditional entropy but not its expectation.

Another point of view is to measure the quantity of information gained by the observation by the Kullback-Leibler risk of the prior relative to the posterior: ${\rm IG}(\bar \Ob_n\rightarrow \theta)=\KL(\pi(\theta)|p(\theta|\bar \Ob_n))$. This quantity is always positive or null and is thus a better candidate to measure a gain in information. Here we can see the difference between information and uncertainty (as in the plane crash example): there is always a gain of information but this  may result in an increase of uncertainty. However when taking expectations, the two points of view meet: the expected information gain is equal to the expected uncertainty reduction.

All this assumes that the model is well specified.

\subsection{Measuring information on the predictive distribution}

The predictive distribution is:
$$ f^{\bon}(x)=\int f^{\theta}(x)p(\theta|\bon) \dd \theta.$$
We wish that $f^{\bon}$ be as close as possible to the true distribution $f^*$. It is natural to consider the Kullback-Leibler risk of $f^{\bon}$ relative to $f^*$, and an estimate of this risk can be used for selection of a Bayesian model. In fact the situation is the same as for the selection of frequentist estimators because we can consider  $f^{\bon}$ as an estimator of $f^*$. For Bayesian model choice we have to estimate the Kullback-Leibler risk, or more easily, the cross-entropy of $f^{\bon}$ relative to $f^*$:
$$\CE(f^{\bon};f^*)=\Eet(\log \frac{1}{f^{\bon}}).$$

This is the quantity that \citet{watanabe2010asymptotic} seeks to estimate by the widely applicable information criterion (WAIC); see also \citet{vehtari2012survey}.
A natural estimator of $\CE(f^{\bon};f^*)$ can be constructed by cross-validation:
$$\CVCE(f^{\Ob})=-\nn \sum_{i=1}^n \log f^{\bar \Ob_{-i}}(X_i),$$
where $\bar \Ob_{-i}$ is the sample from which observation $i$ has been excluded, and $f^{\Ob}$ stands for the set of posterior probabilities.

In fact, this criterion is closely related to the pseudo-Bayes factor. As described in \citet{lesaffre2012bayesian}, the pseudo-Bayes factor for comparing two models $M_1$ and $M_2$ characterized by the sets of  posterior densities $f^{\Ob}_1$ and $f^{\Ob}_2$ respectively, is:
$${\rm PSBF}_{12}=\frac{\prod_{i=1}^n f_1^{\bar \Ob_{-i}}(X_i)}{\prod _{i=1}^n f_2^{\bar \Ob_{-i}}(X_i)}.$$
The pseudo-Bayes factor was introduced by \citet{geisser1979predictive}, and \citet{geisser1980discussion} called the cross-validated density $f_1^{\bar \Ob_{-i}}(X_i)$ the ``conditional predictive ordinate'' (CPO$_i$). We have
$$\log {\rm PSBF}_{12}= \sum_{i=1}^n \log f_1^{\bar \Ob_{-i}}(X_i)-\sum_{i=1}^n \log f_2^{\bar \Ob_{-i}}(X_i)=n[\CVCE(f_1)-\CVCE(f_2)].$$
Thus, the choice between two models using CVCE is equivalent with the choice using pseudo-Bayes factors. \citet{gelfand1994bayesian} showed that PSBF was asymptotically related to AIC; this is not surprising in view of the fact that asymptotically the influence of the prior disappears, and that both criteria are based on estimating the cross-entropy of a predictive distribution with respect to the true distribution.

At first sight CVCE seems to be computationally demanding.
However, this criterion can be computed rather easily using the trick of developing $[f^{\bar \Ob_{-i}}(X_i)]^{-1}$, as shown in \citet{lesaffre2012bayesian} (page 293).
We obtain:
\begin{equation}\CVCE=\nn \si \log \int \frac{1}{f^{\theta}(X_i)}p(\theta|\bar \Ob_n) \dd \theta,\end{equation}
where the integrals can be computed by MCMC. A large number $K$ of realization $\theta_k$ of the parameters values from the posterior distribution can be generated, so that CVCE can be computed as:
$$\CVCE(f^{\Ob})= \nn \sum_{i=1}^n \log K^{-1} \sum_{k=1}^K \frac{1}{f^{\theta_k}(X_i)}.$$
So CVCE, as pseudo-Bayes factors, can be computed with arbitrary precision for any Bayesian model.
 A similar cross-validated criterion was proposed in a prognosis framework by \citet{rizopoulos2015personalized}.

\section{Conclusion}
It is fascinating that the same function, entropy, plays a fundamental role in three different domains, physics, communication theory and statistical inference. Although it is true that these three domains use probability theory as a basis, they remain different and the interpretation of entropy is also different, being a measure of disorder, information and uncertainty respectively. These three concepts are linked but are different. Uncertainty is not just the opposite of information; in some cases uncertainty can increase with additional information.

In statistics, the Kullback-Leibler risk and the cross-entropy are the most useful quantities. The fact that the Kullback-Leibler risk is not a distance is well adapted to the statistical problem: the risk is computed with respect to the true probability, and there is no symmetry between the true probability and a putative probability.
Both estimation and model selection can be done by minimizing a simple estimator of cross-entropy. Among advantages of this point of view is the possibility to interpret a difference of normalized AIC as an estimator of a difference of Kullback-Leibler risks \citep{commenges2008estimating}.
Extensions of the Kullback-Leibler risk by defining the densities on different sigma-fields allows to tackle  several non-standard problems.

Finally, information theory is also relevant in Bayesian theory. Particularly, model choice based on a crossvalidated estimator of the cross-entropy of the predictive distribution is shown to be equivalent to the pseudo-Bayes factor.  \vspace{2ex}

{\bf Acknowledgement} I thank Dimitris Rizopoulos who directed my attention toward pseudo-Bayes factors and the way they can be computed.
\bibliographystyle{Chicago}
\bibliography{Information}

\end{document}